\documentclass[review]{elsarticle}

\usepackage{lineno,hyperref}
\modulolinenumbers[5]

\journal{Statistics and Probability Letter}

\newcommand{\vect}{{\mbox{vec}}}

\newtheorem{theorem}{Theorem}[section]

\usepackage{amssymb,epsfig,verbatim,amsmath}
\usepackage[usenames, dvipsnames]{color}
\usepackage{lscape}
\usepackage{pdflscape}
\usepackage{afterpage}
\usepackage{capt-of}

\begin{document}

\begin{frontmatter}

\title{Semiparametric Estimation for Cure Survival Model with Left-Truncated and Right-Censored Data and Covariate Measurement Error}

\author{Li-Pang Chen\fnref{myfootnote}}
\address{Department of Statistics and Actuarial Science,  University of Waterloo \\
200 University Ave W, Waterloo, ON N2L 3G1}
\fntext[myfootnote]{Corresponding author}
\ead{L358CHEN@uwaterloo.ca}

\begin{abstract}
In this paper, we mainly discuss the cure model with survival data. Different from the usual survival data with right-censoring, we incorporate the features of left-truncation and measurement error in covariates. Generally speaking, left-truncation causes a biased sample in survival analysis; measurement error in covariates may incur a tremendious bias if we do not deal with it properly. To deal with these challenges, we propose a flexible way to analyze left-truncated survival data and correct measurement error in covariates. The theoretical results are also established in this paper.
\end{abstract}


\end{frontmatter}


\section{Introduction} \label{Introduction}

\subsection{Literature Review}
In this paper, our main interest is survival data with cure. In this dataset, these exists a group of subjects who are cured and never experience the failure event (death) in the study period. In the early discussion with right-censored survival data, \cite{Lu2004} considered the semiparametric model.

In the recent developments of cure model, \textit{left-truncation} and \textit{measurement error} are two important features which attract our attention. Left-truncation makes a biased sample in survival data, and measurement error incurs a tremendous bias of the estimator if it is ignored. It is undoubted that these two features make the analysis be challenging.

In the past literature, Chen et al. (2017) proposed the conditional likelihood function based on left-truncation but without measurement error in covariates. With the absence of left-truncation, Ma and Yin (2008) considered the Cox model and introduced a corrected score approach to deal with measurement error in the covariates, but their method can only deal with the linear term of the covariate. To give a more flexible method, Bertrand et al. (2017) implemented the simulation-extrapolation (SIMEX) method which can be used for any function of the covariates.

In many practical situations, these two features may appear in the dataset simultaneously and it may cause the analysis to become complicated and challenging. To the best of our knowledge, there is no method to analyze survival data with those two features incorporated. In this paper, we mainly explore this important problem. We consider the transformation model which includes the Cox model as a special case.


\subsection{Notation and Models}

Let $\xi$ be the calendar time of the recruitment and let $u$ and $r$ denote the calendar time of the initiating event (or the disease incidence) and the failure event, respectively, where $u<r$, and $u < \xi < r$.  Then for those uncored subjects, let $T^{*} = r-u$ be the failure time, and let $ A^{*} = \xi -u$ denote the truncation time. Let $C$ denote the residual censoring time which is measured from $\xi$ to censoring. With both cured and uncored subjects, the failure time is determined by $\widetilde{T} = \pi^\ast T^\ast + (1-\pi^\ast) \infty$, where $\pi^\ast \in \{0,1\}$ indicates whether a subject is cured $(\pi^\ast = 0)$ or not $(\pi^\ast = 1)$. To characterize $\pi^\ast$, we consider a logistic regression model
\begin{eqnarray} \label{Logit-Model}
P\left(\pi^\ast = 1 | Z^\ast \right) = \frac{\exp\left(Z^\ast{}^\top \gamma \right)}{1+\exp\left(Z^\ast{}^\top \gamma \right)},
\end{eqnarray}
where $Z^\ast$ is a $q$-dimensional vector of covariates associated with model (\ref{Logit-Model}), and $\gamma$ is a $q$-dimensional vector of parameters. For subjects who are not cured, we consider the transformation model, which is given by
\begin{eqnarray} \label{Trans-Model}
H\left( T^\ast\right) = -{X^\ast}^\top \beta + \epsilon,
\end{eqnarray}
where $H(\cdot)$ is an unknown increasing function, $\epsilon$ is a random variable with a known distribution, $X^\ast$ is a $p$-dimensional vector of covariates, and $\beta$ is a $p$-dimensional vector of parameters. Model (\ref{Trans-Model}) gives a broad class of some frequently used models in survival analysis. Specifically, when $\epsilon$ has an extreme value distribution, then $T^\ast$ follows the proportional hazards (PH) model; whereas when $\epsilon$ has a logistic distribution, then $T^\ast$ follows the proportional odds (PO) model.

Let $(T,A,X,Z,\pi)$ denote the observed failure time, truncated time, and two covariates which satisfy $\widetilde{T} > A^\ast$. That is, $(T,A,X,Z,\pi) \equiv (\widetilde{T},A^\ast,X^\ast,Z^\ast,\pi^\ast) | \widetilde{T} \geq A^\ast$. For a recruited subject, define $Y = \min\{T,A+C \}$ and $\delta = I\left( T \leq A+C \right)$.

In practice, the covariate $X$ can not be measured correctly and instead we only have an observed covariate $W$. To characteristic the relationship between $X$ and $W$, the classical linear measurement error is frequently used, which is given by
\begin{eqnarray} \label{meas-model}
W = X + \eta,
\end{eqnarray}
where $\eta$ follows the normal distribution with mean zero and covariance matrix $\Sigma_\eta$, and  is independent to $X$. If $\Sigma_\eta$ is unknown, then it can be estimated by additional information, such as repeated measurement or validation data (e.g., \cite{Carroll2006}).
To focus on presenting our proposed method and easing the discussion, we assume that $\Sigma_\eta$ is known.

\subsection{Organization of This Paper}
The remainder is organized as follows. In Section~\ref{Main-Result}, we first present the proposed method to correct the error effect and derive the estimator. After that, we develop the theoretical result for the proposed method. Numerical results are provided in Section~\ref{Numerical}. Finally, we conclude the paper with discussions in Section~\ref{Summary}.


\section{Main Results} \label{Main-Result}

\subsection{Corrected Estimating Equations} \label{Method-SIMEX}

 Suppose that we have an observed sample of $n$ subjects where for $i=1, \cdots, n$, $(Y_{i}, A_{i}, \delta_{i}, W_{i}, Z_i,\pi_i)$ has the same distribution as $\left( Y,A,\delta,W ,Z,\pi\right)$. Let $R_i(t) = I\left( A_i \leq t \leq Y_i \right)$ and $N_i(t) = I \left( Y_i \leq t, \delta_i = 1 \right)$ for $i = 1,\cdots,n$.

As presented in Section~\ref{Introduction}, the covariates $X_i$ is usually unobservable, and instead, we only observe $W_i$. To deal with the mismeasurement and reduce the bias of the estimator, we propose the simulation-extrapolation (SIMEX) method (e.g., \cite{Cook1994}). The proposed procedure is in the following three stages:

\begin{description}
\item[Stage 1]Simulation\\
Let $B$ be a given positive integer and let $\mathcal{Z} = \left\{ \zeta_0,\zeta_1,\cdots,\zeta_M \right\}$ be a sequence of pre-specified values with $0 = \zeta_0<\zeta_1<\cdots<\zeta_M$. where $M$ is a positive integer, and $\zeta_M$ is pre-specified positive number such as $\zeta_M = 1$.

For a given subject $i$ with $i=1,\cdots,n$ and $b = 1,\cdots,B$, we generate $\eta_{b,i}$ from $N\left(0,\Sigma_\eta\right)$. Then for observed vector of covariates $W_i$, we define $W_i(b,\zeta)$ as
\begin{eqnarray} \label{SIMEX-1}
W_i(b,\zeta) = W_i + \sqrt{\zeta} \eta_{b,i}
\end{eqnarray}
for every $\zeta \in \mathcal{Z}$. Therefore, the conditional distribution of $W_i(b,\zeta)$ given $X_i$ is $N\left( X_i, (1+\zeta) \Sigma_\eta \right)$.

\item[Stage 2] Estimation\\
By the similar derivations in \cite{Lu2004}, under left-truncated survival data, we have
\begin{eqnarray} \label{Sec-2-1-1}
&\ \ & P\left( T_i > t | Z_i,W_i(b,\zeta) \right) \nonumber \\
&=& \frac{1}{1+\exp(Z_i^\top \gamma)} + \frac{\exp(Z_i^\top \gamma)}{1+\exp(Z_i^\top \gamma)} \exp\left[ -\Lambda_\epsilon \left\{H(t) + W_i^\top(b,\zeta) \beta \right\} \right] \nonumber \\
&=& \frac{\bar{G}(Z_i^\top \gamma)}{G(\Lambda_\epsilon \left\{H(t) + W_i^\top(b,\zeta) \beta \right\} - Z_i^\top \gamma)},
\end{eqnarray}
where $\Lambda_\epsilon(\cdot)$ is the cumulative hazard function of $\epsilon$, $G(x) = \frac{\exp(x)}{1+\exp(x)}$, and $\bar{G}(x) = 1-G(x)$. Taking log function with negative sign on (\ref{Sec-2-1-1}) gives
\begin{eqnarray} \label{Sec-2-1-2}
&\ \ & -\log\left\{ P\left( T_i > t | Z_i,W_i(b,\zeta) \right) \right\} \nonumber \\
&=& \log \left\{ G(\Lambda_\epsilon \left\{H(t) + W_i^\top(b,\zeta) \beta \right\} - Z_i^\top \gamma) \right\} - \log \left\{ \bar{G}(Z_i^\top \gamma) \right\}.
\end{eqnarray}

By the counting process techniques (e.g., \cite{Anderson1993}), we define
\begin{eqnarray} \label{Sec-2-1-3}
M_i(t) = N_i(t) - \int_0^t R_i(u) d \log \left\{  G(\Lambda_\epsilon \left\{H(u) + W_i^\top(b,\zeta) \beta \right\} - Z_i^\top \gamma) \right\},
\end{eqnarray}
which is a martingale process with $E\left\{M_i(t) \right\} = 0$. Then based on (\ref{Sec-2-1-3}), we have two estimating equations (EE):
\begin{eqnarray} \label{EE-SIMEX-1} 
\sum \limits_{i=1}^n \left[ dN_i(t) -  R_i(t) d \log \left\{  G(\Lambda_\epsilon \left\{H(t) + W_i(b,\zeta)^\top \beta \right\} - Z_i^\top \gamma) \right\} \right] = 0
\end{eqnarray}
and
\begin{eqnarray} \label{EE-SIMEX-2-pre}
\sum \limits_{i=1}^n \int_0^\infty W_i(b,\zeta) \left[ dN_i(t) -  R_i(u) d \log \left\{  G(\Lambda_\epsilon \left\{H(t) + W_i(b,\zeta)^\top \beta \right\} - Z_i^\top \gamma) \right\} \right] = 0.
\end{eqnarray}
Let $\theta = (\beta^\top, \gamma^\top)^\top$ be a $(p+q)$-dimensional vector of parameters. Solving (\ref{EE-SIMEX-1}) yields the estimator of $H(\cdot)$ when both $\beta$ and $\gamma$ are fixed, which is denoted by $\widehat{H}(t;b,\zeta,\theta)$. However, (\ref{EE-SIMEX-2-pre}) only gives the estimator of $\beta$. To derive the estimator of $\gamma$, we need to develop the third estimating equation based on $\pi_i$.
We consider the conditional probability
\begin{eqnarray*}
&\ \ &P\left( \pi_i = 1 | \delta_i, T_i, W_i(b,\zeta), Z_i \right) \\
&=& \delta_i + (1-\delta_i) \bar{G}(\Lambda_\epsilon \left\{H(t) + W_i^\top(b,\zeta) \beta \right\} - Z_i^\top \gamma),
\end{eqnarray*}
and by the similar derivation of Equation (9) in \cite{Lu2004}, we have the unbiased estimating equation for $\gamma$:
\begin{eqnarray}  \label{EE-SIMEX-3-pre}
\sum \limits_{i=1}^n W_i(b,\zeta) \left(\delta_i + (1-\delta_i) \bar{G}\left[\Lambda_\epsilon \left\{H(T_i) + W_i^\top(b,\zeta) \beta \right] - G(Z_i^\top \gamma)\right\} \right).
\end{eqnarray}

Replacing $H(\cdot)$ in (\ref{EE-SIMEX-2-pre}) and (\ref{EE-SIMEX-3-pre}) by $\widehat{H}(t;b,\zeta,\theta)$ gives the following two estimating equations 
\begin{eqnarray} \label{EE-SIMEX}
&\ \ &U_{SIMEX,1}(\theta) \nonumber \\
&\triangleq& \frac{1}{n} \sum \limits_{i=1}^n \int_0^\infty W_i(b,\zeta) \left[ dN_i(t) \right. \nonumber \\
&\ \ &-\left.  R_i(u) d \log \left\{  G(\Lambda_\epsilon \left\{\widehat{H}(t;b,\zeta,\theta) + W_i(b,\zeta)^\top \beta \right\} - Z_i^\top \gamma) \right\} \right] \nonumber  \\
&=& 0
\end{eqnarray}
and
\begin{eqnarray} \label{EE-SIMEX-2}
&\ \ &U_{SIMEX,2}(\theta) \nonumber \\
&\triangleq& \frac{1}{n} \sum \limits_{i=1}^n W_i(b,\zeta) \left(\delta_i + (1-\delta_i) \bar{G}\left[\Lambda_\epsilon \left\{\widehat{H}(T_i;b,\zeta,\theta) + W_i(b,\zeta)^\top \beta \right\} - G(Z_i^\top \gamma)\right] \right) \nonumber  \\
&=& 0
\end{eqnarray}
for every $b = 1,\cdots,B$ and $\zeta \in \mathcal{Z}$. Let $\widehat{\theta}(b,\zeta)$ denote the solution of two estimating equations $U_{SIMEX,1}(\theta) = 0$ and $U_{SIMEX,2}(\theta) = 0$. Moreover, we define
\begin{eqnarray} \label{SIMEX-2}
\widehat{\theta}(\zeta) = \frac{1}{B} \sum \limits_{b=1}^B \widehat{\theta}(b,\zeta).
\end{eqnarray}

\item[Stage 3] Extrapolation\\
By (\ref{SIMEX-2}), we have a sequence $\left\{ \left( \zeta, \widehat{\theta}(\zeta)\right) : \zeta \in \mathcal{Z} \right\}$. Then we fit a regression model to the sequence
\begin{eqnarray}
\widehat{\theta}(\mathcal{Z}) = \varphi\left(\mathcal{Z}, \Gamma\right) + \varrho,
\end{eqnarray}
where $\varphi(\cdot)$ is the user-specific regression function, $\Gamma$ is the associated parameter, and $\varrho$ is the noise term. The parameter $\Gamma$ can be estimated by the least square method, and we let $\widehat{\Gamma}$ denote the resulting estimate of $\Gamma$.

Finally, we calculate the predicted value
\begin{eqnarray} \label{SIMEX-3}
\widehat{\theta}_{SIMEX} = \varphi \left(-1, \widehat{\Gamma} \right)
\end{eqnarray}
and take $\widehat{\theta}_{SIMEX}$ as the \textit{SIMEX} estimator of $\theta$.

\item[Stage 4] Estimation of $H(\cdot)$ \\
Furthermore, we can also derive the estimator of the unknown function $H(\cdot)$. To do this, we first replace $\theta$ by $\widehat{\theta}_{SIMEX}$ in $\widehat{H}(t;b,\zeta,\theta)$, which gives $\widehat{H}(t;b,\zeta,\widehat{\theta}_{SIMEX})$. For every $t$ and $\zeta$, taking average with respect to $b$ gives $\widehat{H}(t;\zeta,\widehat{\theta}_{SIMEX}) = \frac{1}{B} \sum \limits_{b=1}^B \widehat{H}(t;b,\zeta,\widehat{\theta}_{SIMEX})$. Finally, similar to Stage 3 above, fitting a regression model and taking $\zeta = -1$ as a predicted value yields a final estimator $\widehat{H}(t;\widehat{\theta}_{SIMEX})$, also denoted as $\widehat{H}_{SIMEX}(t)$.
\end{description}

\subsection{Theoretical Results}

In this section, we present the theoretical results of the proposed method. We first define some notation. Let $\theta_0 = \left( \beta_0^\top, \gamma_0^\top \right)^\top$ denote the true value of the parameter $\theta$, and let $H_0(\cdot)$ denote the true function of $H(\cdot)$. Let $\lambda_\epsilon(t) = \frac{d \Lambda_\epsilon(t)}{dt}$. For $i=1,\cdots,n$, define
\begin{eqnarray*}
\zeta_i(t;X_i,H,\theta) &=& \frac{d}{dt} \log \left[ \lambda_\epsilon \left\{ H(t) + X_i^\top \beta \right\} \right] \\
&\ \ & - \lambda_\epsilon \left\{ H(t) + X_i^\top \beta \right\} G(\Lambda_\epsilon \left\{H(t) + X_i^\top \beta \right\} - Z_i^\top \gamma), \\
\Phi_{1i}(t;X_i,H,\theta) &=& \left(X_i^\top, Z_i^\top G(\Lambda_\epsilon \left\{H(t) + X_i^\top \beta \right\} - Z_i^\top \gamma) \right)^\top, \\
\Phi_{2i}(t;X_i,H,\theta) &=& \left(X_i^\top \zeta_i(t;X_i,H,\theta), Z_i^\top G(\Lambda_\epsilon \left\{H(t) + X_i^\top \beta \right\} - Z_i^\top \gamma) \right)^\top, \\
\Phi_{3i}(t;X_i,H,\theta) &=& \left(X_i^\top, Z_i^\top (1-\delta_i) G(\Lambda_\epsilon \left\{H(t) + X_i^\top \beta \right\} - Z_i^\top \gamma) \right)^\top, \\
\Psi_i(t;X_i,H,\theta) &=& \lambda_\epsilon \left\{ H(t) + X_i^\top \beta \right\} \bar{G}(\Lambda_\epsilon \left\{H(t) + X_i^\top \beta \right\} - Z_i^\top \gamma).
\end{eqnarray*}
We further define
\begin{eqnarray*}
B(t,s;X_i) &=& \exp \left( \int_s^t \frac{E\left\{ \zeta_i (u;X_i,H_0,\theta_0) dN_i(u) \right\}}{E\left\{ \Psi_i(u;X_i,H_0,\theta_0) R_i(u) \right\}} \right) \\
\phi_i(t;X_i) &=& \frac{E\left\{ \Phi_{3i}(Y_i;X_i,H_0,\theta_0) \Psi_i(Y_i;X_i,H_0,\theta_0) R_i(t) B(t,Y_i;X_i) \right\}}{E\left\{ \Psi_i(t;X_i,H_0,\theta_0) R_i(t) \right\}}.
\end{eqnarray*}

 We now present the theoretical results of $\widehat{\theta}_{SIMEX}$ and $\widehat{H}_{SIMEX}(\cdot)$ in the following theorem.
\begin{theorem} \label{Thm-SIMEX}
Under regularity conditions in \ref{R.C.}, estimators $\widehat{\theta}_{SIMEX}$ and $\widehat{H}_{SIMEX}(\cdot)$ have the following properties:
\begin{itemize}
\item[(1)] $\widehat{\theta}_{SIMEX} \stackrel{p}{\longrightarrow} \beta_0$  as  $n \rightarrow \infty$;
\item[(2)] $\widehat{H}_{SIMEX}(t) \stackrel{p}{\longrightarrow} H_0(t)$  as  $n \rightarrow \infty$;
\item[(3)] $\sqrt{n} \left(\widehat{\theta}_{SIMEX} - \beta_0 \right) \stackrel{d}{\longrightarrow} N\left(0, \left\{\frac{\partial \varphi}{\partial \Gamma}\left(-1, \widehat{\Gamma} \right) \right\} \mathcal{Q} \left\{ \frac{\partial \varphi}{\partial \Gamma} \left(-1, \widehat{\Gamma} \right) \right\}^\top  \right)$  as  $n \rightarrow \infty$;
\item[(4)] $\sqrt{n} \left\{ \widehat{H}_{SIMEX}(t) - H_0(t) \right\}$ converges to the Gaussian process with mean zero and covariance function $E\left\{ \mathcal{H}_i(t) \mathcal{H}_i(s) \right\}$,
\end{itemize}
where the exact formulations of $\mathcal{Q}$ and $\mathcal{H}_i(t)$ are placed in \ref{pf}.
\end{theorem}

\section{Numerical Study} \label{Numerical}

\subsection{Simulation Setup}

We examine the setting where $\epsilon$ is generated from the extreme value distribution and the logistic distribution, and the truncation time $A^\ast$ is generated from the exponential distribution with mean one. Let $\theta = (\beta, \gamma)^\top$ denote a two-dimensional vector of parameters, and let $\theta_0 = (\beta_0, \gamma_0)^\top$ be the vector of true parameters where we set $\theta_0 = \left( 1,1 \right)^\top$. We consider a scenario where $ (X^\ast, Z^\ast)^\top$ are generated from a bivariate normal distribution with mean zero and variance-covariance matrix $\Sigma$, which is set as $ \left( \begin{array}{ c c} 
4 & 0.7 \\
0.7 & 3
\end{array}  
   \right)$. 
 Given $\epsilon$, $X^\ast$ and $\beta_0$, the failure time $T^\ast$ is generated from the model: 
\begin{eqnarray*} 
\log T^\ast = -  X^\ast \beta_0  + \epsilon.
\end{eqnarray*}
Based on our two settings of $\epsilon$, the failure time $T^\ast$ follows the PH model and the PO model, respectively. On the other hand, $\pi$ is generated by (\ref{Logit-Model}), and hence, the failure time with cure is determined by $\widetilde{T} = \pi^\ast T^\ast + (1-\pi^\ast) \infty$.  Therefore, the observed data $(A,T,X,Z)$ is collected from $(A^\ast,\widetilde{T},X^\ast,Z^\ast)$ by conditioning on that $\widetilde{T} \geq A^\ast$. We repeatedly generate data these steps we obtain a sample of a required size $n=200$. For the measurement error process, we consider model $(\ref{meas-model})$  with error $\eta \sim N \left( 0, \Sigma_\eta \right)$, where  $\Sigma_\eta$ is a scalar which is taken as $ 0.01$, $0.5$, and $0.75$, respectively.
  
We consider two censoring rates, say 25\% and 50\%, and let the  censoring time $C$ be generated from the uniform distribution $U(0,c)$, where $c$ is determined by a given censoring rate.
Consequently, $Y$ and $\delta$ are determined by $Y = \min \left\{ T, A+C \right\}$ and $\delta = I \left( T \leq A + C \right)$. In implementing the proposed method, we set $B = 500$ and partition the interval $[0,2]$ into subintervals with width $0.25$, and let the resulting cutpoints be the values of $\zeta$. We take the regression function $\varphi(\cdot)$ to be the quadratic polynomial function, which is a widely used function in many cases (e.g., Cook and Stefaski 1994; Carroll et al. 2006). Finally, 1000 simulations are run for each parameter setting.

\subsection{Simulation Results}

We mainly examine the performance of the proposed method which is denoted by Chen ($\widehat{\theta}_{SIMEX}$). In addition, to see the impact of the measurement error in covariate, we examine the \textit{naive estimator} which is obtained by implementing $W_i$ in the estimating equations instead of $W_i(b,\zeta)$, and the naive estimator is denoted by Naive ($\widehat{\theta}_{naive}$).  We report the biases of estimates (Bias), the empirical variances (Var), the mean squared errors (MSE), and the coverage probabilities (CP) of those two estimators. The results are reported in Table~\ref{tab:Sim}. 

First, the censoring rate and measurement degree have noticeable impact on each estimation methods.
As expected, biases and variance estimates increase as the censoring rate increases. 
When the measurement degree increases, biases of both $\widehat{\theta}_{naive}$ and $\widehat{\theta}_{SIMEX}$ are increasing, and the impact of the measurement error degrees seems more obvious on the naive estimator $\widehat{\theta}_{naive}$.

Within a setting with a given censoring rate and a measurement error degree, the naive method and the proposed method perform differently. When measurement error occurs, the performance of the proposed method is better than the naive method. The naive method produces considerable finite sample biases with coverage rates of 95\% confidence intervals significantly departing from the nominal level. The proposed method outputs satisfactory estimate with small finite sample biases and reasonable coverage rates of 95\% confidence intervals. Compared to the variance estimates produced by the naive approach, the proposed method which accounts for measurement error effects yield larger variance estimates, and this is the price paid to remove biases in point estimators. This phenomenon is typical in the literature of measurement error models. However, mean squared errors produced by the proposed method tends to be a lot smaller than those obtained from the naive method. 
 
\section{Discussion} \label{Summary}

In this article, we focus the discussion on the transformation model based on cured survival data with left-truncation and develop a valid method to correct the covariate measurement error and derive an efficient estimator. In this article, we also establish the large sample properties, and the numerical results guarantee that our proposed method outperforms. Although we only focus on the simple structure of the measurement error model and assume that $Z_i$ is precisely measured, our method can easily be extended to complex measurement error models or additional information, such as repeated measurement or validation data, and also allows $Z_i$ in (\ref{Logit-Model}) is mismeasured. In addition, there are still many challenges in this topic, such as the discussion of time-dependent covariates with mismeasurement. These topics are also our researches in the future.

 \begin{table}
       \huge
     \caption{Numerical results for simulation study} \label{tab:Sim}

 \scriptsize

\center
   \renewcommand{\arraystretch}{0.75}
 \begin{tabular}{c c c c  ccccccccccccc}

 \\
 \hline\hline
model & cr  & $\sigma_{\eta}$  &  Method  & \multicolumn{4}{c} {Estimator of $\beta$ } & & \multicolumn{4}{c}{Estimator of $\gamma$ } \\ \cline{5-8}  \cline{10-13}

 &  &  &  & Bias & Var & MSE & CP(\%) &  & Bias & Var  & MSE & CP(\%)
\\
 \hline 
PH & 25\% &  0.01 & Naive$(\widehat{\beta}_{naive})$ & -0.230 & 0.007 & 0.059 & 21.3 &  & -0.749 & 0.014 &  0.626 & 14.9\\
  &    &     & Chen$(\widehat{\beta}_{SIMEX})$    & 0.017 & 0.013  & 0.014 & 94.7 &  & 0.009 & 0.028 &  0.028 & 94.2   \\
      \\
  &    & 0.50 & Naive$(\widehat{\beta}_{naive})$ & -0.343 & 0.006  & 0.123 & 1.6 &  & -0.606 & 0.015 & 0.432 & 30.0\\
  &    &       & Chen$(\widehat{\beta}_{SIMEX})$        & 0.025 & 0.023 & 0.026 & 94.5 &  & 0.011 & 0.027  & 0.028 & 94.5\\
      \\
  &    & 0.75  & Naive$(\widehat{\beta}_{naive})$ &-0.347 & 0.005  & 0.125 & 0.3 &  & -0.636 & 0.016  & 0.465 & 23.8\\
  &    &       & Chen$(\widehat{\beta}_{SIMEX})$       & 0.025 & 0.023  & 0.023 & 94.8 &  & 0.019 & 0.025  & 0.025 & 93.9\\       
& 50\% & 0.01 & Naive$(\widehat{\beta}_{naive})$ & -0.248 & 0.016  & 0.267 & 9.1 &  & -0.742 & 0.016  & 0.565 & 0.1\\
 &     &     & Chen$(\widehat{\beta}_{SIMEX})$       & 0.017 & 0.014  & 0.014 & 94.4 &  & 0.016 & 0.021  & 0.021 & 94.3\\
      \\
  &    & 0.50 & Naive$(\widehat{\beta}_{naive})$ &  -0.375 & 0.015  & 0.145 & 0.2 &  & -0.600 & 0.016  & 0.376 & 0.4\\
   &   &       & Chen$(\widehat{\beta}_{SIMEX})$        & 0.024 & 0.036  & 0.039 & 95.2 &  & 0.019 & 0.025 & 0.025 & 95.0\\
      \\
 &      & 0.75  & Naive$(\widehat{\beta}_{naive})$ & -0.360 & 0.014  & 0.134 & 0.1 &  & -0.630 & 0.014  & 0.413 & 0.2\\
 &     &       & Chen$(\widehat{\beta}_{SIMEX})$        &0.025 & 0.033 & 0.033 & 94.6 &  & 0.026 & 0.025  & 0.025 & 94.8\\   
  \hline 
PO & 25\% &  0.01 & Naive$(\widehat{\beta}_{naive})$ &  -0.250 & 0.009 &  0.072 & 23.0 & &  -0.729 & 0.015  & 0.557 & 0.4 \\
  &    &     & Chen$(\widehat{\beta}_{SIMEX})$    & 0.010 & 0.019  & 0.020 & 94.2 &  & 0.009 & 0.024  & 0.024 & 94.5   \\
      \\
  &    & 0.50 & Naive$(\widehat{\beta}_{naive})$ & -0.377 & 0.008  & 0.150 & 1.5 &  & -0.588 & 0.017  & 0.369 & 3.6\\
  &    &       & Chen$(\widehat{\beta}_{SIMEX})$        & 0.012 & 0.018  & 0.040 & 94.5 &  & 0.011 & 0.024 & 0.025 & 93.7\\
      \\
  &    & 0.75  & Naive$(\widehat{\beta}_{naive})$ &-0.362 & 0.007  & 0.138 & 1.1 &  & -0.619 & 0.014  & 0.405 & 1.4\\
  &    &       & Chen$(\widehat{\beta}_{SIMEX})$       & 0.016 & 0.018  & 0.018 & 94.3 &  & 0.015 & 0.022  & 0.022 & 94.6\\       
& 50\% & 0.01 & Naive$(\widehat{\beta}_{naive})$ & -0.268 & 0.016  & 0.273 & 10.1 &  & -0.842 & 0.016  & 0.574 & 1.3\\
 &     &     & Chen$(\widehat{\beta}_{SIMEX})$       & 0.016 & 0.024  & 0.024 & 94.6 &  & 0.016 & 0.027  & 0.027 & 94.5\\
      \\
  &    & 0.50 & Naive$(\widehat{\beta}_{naive})$ &  -0.388 & 0.016  & 0.168 & 1.4 &  & -0.600 & 0.016  & 0.376 & 1.4\\
   &   &       & Chen$(\widehat{\beta}_{SIMEX})$        & 0.027 & 0.036  & 0.037 & 94.2 &  & 0.021 & 0.026 & 0.026 & 95.1\\
      \\
 &      & 0.75  & Naive$(\widehat{\beta}_{naive})$ & -0.410 & 0.017  & 0.185 & 1.9 &  & -0.630 & 0.018  & 0.413 & 1.2\\
 &     &       & Chen$(\widehat{\beta}_{SIMEX})$        &0.028 & 0.036 & 0.036 & 94.6 &  & 0.025 & 0.027  & 0.027 & 94.6\\   

 \hline\hline
 \\
\end{tabular}
\\
\tiny
\raggedright

Note:
\\ $\times$ - usage of the true covariate $X$;
\\ cr - censoring rate;
\\Bias - Difference between empirical mean and true value;
\\Var - Empirical variance;
\\MSE - Mean square error;
\\MVE - Model-based variance;
\\CP - Model-based coverage probability.
\end{table}

\clearpage





\appendix

\section{Regularity Conditions} \label{R.C.}
\begin{itemize}
\item[(C1)] $\Theta$ is a compact set, and the true parameter value $\theta_0$ is an interior point of $\Theta$.
\item[(C2)] Let $\tau$ be the finite maximum support of the failure time.
\item[(C3)] The $\left\{ A_i, Y_i,X_i,Z_i \right\}$ are independent and identically distributed for $i=1,\cdots,n$. 
\item[(C4)] The covariates $X_i$ and $Z_i$ are bounded.
\item[(C5)] Conditional on the covariates $X_i^\ast$ and $Z_i^\ast$, $ T_i^\ast$ is independent of $A_i^\ast$.
\item[(C6)] Censoring time $C_i$ is non-informative. That is, the failure time $T_i$ and the censoring time $C_i$ are independent, given the covariates $\{Z_i, X_i\}$.
\item[(C7)] The regression function $\varphi(\cdot)$ is true, and its first order derivative exists.
\end{itemize}
Condition (C1) is a basic condition that is used to derive the maximizer of the target function. (C2) to (C6) are standard conditions for survival analysis, which allow us to obtain the sum of i.i.d. random variables and hence to derive the asymptotic properties of the estimators. Condition (C7) is a common assumption in SIMEX method.


\setcounter{equation}{0}
\renewcommand\theequation{B.\arabic{equation}}
\section{Proof of Theorem~\ref{Thm-SIMEX}} \label{pf} \ \\
\underline{Proof of Theorem~\ref{Thm-SIMEX} (1)}:\\
Let 
\begin{eqnarray} \label{pf-1-1}
U_{SIMEX}(\theta) = \left( U_{SIMEX,1}^\top(\theta), U_{SIMEX,2}^\top(\theta) \right)^\top,
\end{eqnarray}
and let $\theta(b,\zeta)$ denote a solution of $E\left\{U_{SIMEX}(\theta)\right\} = 0$. Since $\widehat{\theta}(b,\zeta)$ is a solution of $U_{SIMEX}(\theta) = 0$. By the Uniformly Law of Large Numbers (e.g., \cite{van-der1998}), we have that $\frac{1}{n} U_{SIMEX}(\theta)$ converges uniformly to $E\left\{U_{SIMEX}(\theta)\right\}$. Then we have that as $n \rightarrow \infty$,
\begin{eqnarray} \label{pf-1-3}
\widehat{\theta}(b,\zeta) \stackrel{p}{\longrightarrow} \theta(b,\zeta).
\end{eqnarray}
By definition (\ref{SIMEX-2}), taking averaging with respect to $b$ on both sides of (\ref{pf-1-3}) gives that as $n \rightarrow \infty$,
\begin{eqnarray}  \label{pf-1-4}
\widehat{\theta}(\zeta) \stackrel{p}{\longrightarrow} \theta(\zeta)
\end{eqnarray}
for every $\zeta \in \mathcal{Z}$.
By (\ref{pf-1-4}), we can show that as $n \rightarrow \infty$,
\begin{eqnarray} 
\widehat{\Gamma} \stackrel{p}{\longrightarrow} \Gamma.
\end{eqnarray}
Since $\widehat{\beta}_{SIMEX} = \varphi\left(-1,\widehat{\Gamma} \right)$, therefore, by the continuous mapping theorem, we have that as $n \rightarrow \infty$,
\begin{eqnarray} \label{pf-1-5}
\widehat{\beta}_{SIMEX} \stackrel{p}{\longrightarrow} \beta_0.
\end{eqnarray}
\ \\
\underline{Proof of Theorem~\ref{Thm-SIMEX} (2)}:\\
By (\ref{pf-1-5}), we have $\widehat{H}(t;b,\zeta,\widehat{\theta}_{SIMEX}) - \widehat{H}(t;b,\zeta,\theta_0) = o_p(1)$ for every $t \in [0,\tau]$, b, and $\zeta$. Taking average with respect to $b$ gives $\widehat{H}(t;\zeta,\widehat{\theta}_{SIMEX}) - \widehat{H}(t;\zeta,\theta_0) = o_p(1)$. On the other hand, by the Uniformly Law of Large Numbers and similar derivations in \cite{Lu2004} with $\zeta \rightarrow -1$, we have that as $n \rightarrow \infty$, $\widehat{H}(t;\theta_0) - H_0(t) \stackrel{p}{\longrightarrow} 0$ for all $t \in [0,\tau]$. Therefore, we conclude that as $n \rightarrow \infty$, $\widehat{H}(t;,\widehat{\theta}_{SIMEX}) - H_0(t) \stackrel{p}{\longrightarrow} 0$ by the fact that $\widehat{H}(t;-1,\widehat{\theta}_{SIMEX}) - H_0(t) = \widehat{H}(t;-1,\widehat{\theta}_{SIMEX}) - \widehat{H}(t;-1,\theta_0) + \widehat{H}(t;-1,\theta_0) - H_0(t)$.
\\
\\
\underline{Proof of Theorem~\ref{Thm-SIMEX} (3)}:\\
For $b = 1,\cdots,B$ and $\zeta \in \mathcal{Z}$, applying the Taylor series expansion on (\ref{pf-1-1}) around ${\theta}(b,\zeta)$ gives
\begin{eqnarray*}
0 &=& U_{SIMEX}\left( \widehat{\theta}(b,\zeta) \right) \\
&=& U_{SIMEX}\left( {\theta}(b,\zeta) \right) + \frac{\partial U_{SIMEX}\left( {\theta}(b,\zeta) \right)}{\partial \theta} \left\{ \widehat{\theta}(b,\zeta) - {\theta}(b,\zeta) \right\} + o_p\left( \frac{1}{\sqrt{n}}\right),
\end{eqnarray*}
or equivalently,
\begin{eqnarray} \label{pf-3-1}
\sqrt{n} \left\{ \widehat{\theta}(b,\zeta) - {\theta}(b,\zeta) \right\} &=& \left(-\frac{\partial U_{SIMEX}\left( {\theta}(b,\zeta) \right)}{\partial \theta} \right)^{-1} \sqrt{n} U_{SIMEX}\left( {\theta}(b,\zeta) \right) \nonumber \\
&\ \ & + o_p\left( 1 \right).
\end{eqnarray}

By (\ref{EE-SIMEX}), (\ref{EE-SIMEX-2}), and the Uniformly Law of Large Numbers, we have that as $n \rightarrow \infty$,
\begin{eqnarray} \label{pf-3-2}
\left(-\frac{\partial U_{SIMEX}\left( {\theta}(b,\zeta) \right)}{\partial \theta} \right) \stackrel{p}{\longrightarrow} \mathcal{A}\left(b,\zeta \right),
\end{eqnarray}
where 
\begin{eqnarray*}
&& \mathcal{A}\left(b,\zeta \right) \\
&=& - E \left[ \int_0^\infty \left\{ \Phi_{1i}(t;W_i(b,\zeta),H_0,\theta_0) - \phi_i(t;W_i(b,\zeta)) \right\} \Phi_{2i}^\top(t;W_i(b,\zeta),H_0,\theta_0) dN_i(t) \right]
\end{eqnarray*}

On the other hand, by (\ref{Sec-2-1-3}), the  estimating equations (\ref{EE-SIMEX}) and (\ref{EE-SIMEX-2}) can be expressed as
\begin{eqnarray*}
U_{SIMEX,1}({\theta}(b,\zeta)) = \frac{1}{n} \sum \limits_{i=1}^n \int_0^\infty \left\{ W_i(b,\zeta) - \phi_{\beta;i}(t;W_i(b,\zeta)) \right\} dM_i(t) + o_p\left( \frac{1}{\sqrt{n}} \right)
\end{eqnarray*}
and 
\begin{eqnarray*}
U_{SIMEX,2}({\theta}(b,\zeta)) &=& \frac{1}{n} \sum \limits_{i=1}^n \int_0^\infty \left\{ Z_i^\top G(\Lambda_\epsilon \left\{H(t) + W_i(b,\zeta)^\top \beta \right\} - Z_i^\top \gamma) \right. \\
&& -\left. \phi_{\gamma;i}(t;W_i(b,\zeta)) \right\} dM_i(t) + o_p\left( \frac{1}{\sqrt{n}} \right),
\end{eqnarray*}
where $\phi_{\beta;i}(\cdot)$ is the first $p$-dimensional components of $\phi_{i}(\cdot)$  and $\phi_{\gamma;i}(\cdot)$ is the remaining $q$-dimensional components of $\phi_{i}(\cdot)$. Thus $U_{SIMEX}({\theta}(b,\zeta))$ can be derived as a sum of i.i.d. random functions, which is given by
\begin{eqnarray}  \label{pf-3-3}
\sqrt{n} U_{SIMEX}({\theta}(b,\zeta)) = \frac{1}{\sqrt{n}} \sum \limits_{i=1}^n \mathcal{B}_i(b,\zeta) + o_p(1),
\end{eqnarray}
where
\begin{eqnarray*}
\mathcal{B}_i(b,\zeta) &=& \int_0^\infty \left\{ \Phi_{1i}(t;W_i(b,\zeta),H_0,\theta(b,\zeta)) - \phi_i(t;W_i(b,\zeta)) \right\} dM_i(t).
\end{eqnarray*}

Combining (\ref{pf-3-3}) and (\ref{pf-3-2}) with (\ref{pf-3-1}) yields
\begin{eqnarray} \label{pf-3-4}
\sqrt{n} \left\{ \widehat{\theta}(b,\zeta) - {\theta}(b,\zeta) \right\} = \frac{1}{\sqrt{n}} \sum \limits_{i=1}^n \mathcal{A}^{-1}\left(b,\zeta \right) \mathcal{B}_i(b,\zeta) + o_p\left( 1 \right).
\end{eqnarray}
By (\ref{SIMEX-2}), taking average with respect to $b$ on both sides of (\ref{pf-3-4}) gives
\begin{eqnarray} \label{pf-3-5}
\sqrt{n} \left\{ \widehat{\theta}(\zeta) - {\theta}(\zeta) \right\} = \frac{1}{\sqrt{n}} \sum \limits_{i=1}^n \mathbf{B}_i(\zeta) + o_p\left( 1 \right)
\end{eqnarray}
for $\zeta \in \mathcal{Z}$, where $ \mathbf{B}_i(\zeta) =  \frac{1}{B} \sum \limits_{b=1}^B \mathcal{A}^{-1}\left(b,\zeta \right) \mathcal{B}_i(b,\zeta)$.

Let $\widehat{\theta}(\mathcal{Z}) = \vect\left\{\widehat{\theta}(\zeta) : \zeta \in \mathcal{Z}\right\}$ denote the vectorization of estimator $\widehat{\theta}(\zeta)$ with every $\zeta \in \mathcal{Z}$. By the Central Limit Theorem on (\ref{pf-3-5}), we have that as $n \rightarrow \infty$,
\begin{eqnarray} \label{pf-3-6}
\sqrt{n} \left\{ \widehat{\theta}(\mathcal{Z}) - {\theta}(\mathcal{Z}) \right\} \stackrel{d}{\longrightarrow} N\left(0, \Omega\left( \mathcal{Z} \right) \right),
\end{eqnarray} 
where $\Omega\left( \mathcal{Z} \right) = \text{cov}\left\{ \mathbf{B}_i(\mathcal{Z}) \right\}$. By the Taylor series expansion on $\varphi\left( \mathcal{Z}, \Gamma \right)$ with respect to $\Gamma$, we have
\begin{eqnarray} \label{pf-3-7}
\varphi\left( \mathcal{Z}, \widehat{\Gamma} \right) - \varphi\left( \mathcal{Z}, \Gamma \right) \approx \frac{\partial \varphi\left( \mathcal{Z}, \Gamma \right)}{\partial \Gamma} \left(\widehat{\Gamma} - \Gamma \right).
\end{eqnarray}
Let $\mathcal{C} = \frac{\partial \varphi\left( \mathcal{Z}, \Gamma \right)}{\partial \Gamma}$ and $\mathcal{D} = \left\{\frac{\partial \varphi\left( \mathcal{Z}, \Gamma \right)}{\partial \Gamma}\right\}^\top \frac{\partial \varphi\left( \mathcal{Z}, \Gamma \right)}{\partial \Gamma}$. Combining (\ref{pf-3-6}) and (\ref{pf-3-7}) gives that as $n \rightarrow \infty$,
\begin{eqnarray} \label{pf-3-8}
\sqrt{n} \left(\widehat{\Gamma} - \Gamma \right) \stackrel{d}{\longrightarrow} N\left(0, \mathcal{D}^{-1} \mathcal{C} \Omega(\mathcal{Z}) \mathcal{C}^\top \mathcal{D}^{-1}  \right).
 \end{eqnarray}
 Finally, since the SIMEX estimator is defined by $\widehat{\beta}_{SIMEX} = \varphi\left(-1, \widehat{\Gamma} \right)$. Let $\mathcal{Q} = \mathcal{D}^{-1} \mathcal{C} \Omega(\mathcal{Z}) \mathcal{C}^\top \mathcal{D}^{-1}$. Combining (\ref{pf-3-7}) and (\ref{pf-3-8}) with $\zeta \rightarrow -1$ and applying the delta method give that as $n \rightarrow \infty$,
 \begin{eqnarray*} 
\sqrt{n} \left(\widehat{\theta}_{SIMEX} - \theta_0 \right) \stackrel{d}{\longrightarrow} N\left(0, \left\{\frac{\partial \varphi}{\partial \Gamma}\left(-1, \widehat{\Gamma} \right) \right\} \mathcal{Q} \left\{ \frac{\partial \varphi}{\partial \Gamma} \left(-1, \widehat{\Gamma} \right) \right\}^\top  \right).
 \end{eqnarray*}
 \ \\ 
 \underline{Proof of Theorem~\ref{Thm-SIMEX} (4)}:\\
 We first consider the expression of $\sqrt{n} \left\{ \widehat{H}(t;b,\zeta,\widehat{\theta}_{SIMEX}) - \widehat{H}(t;b,\zeta,\theta_0) \right\}$. By the Taylor series expansion with respect to $\beta$, we have
 \begin{eqnarray} \label{pf-4-1}
&& \sqrt{n} \left\{ \widehat{H}(t;b,\zeta,\widehat{\theta}_{SIMEX}) - \widehat{H}(t;b,\zeta,\theta_0) \right\} \nonumber \\
&=& \frac{\partial \widehat{H}(t;b,\zeta,\theta_0)}{\partial \theta} \sqrt{n} \left( \widehat{\theta}_{SIMEX} - \theta_0 \right) \nonumber \\ 
&=& A(t) \sqrt{n} \left( \widehat{\theta}_{SIMEX} - \theta_0 \right) + o_p(1) \nonumber \\
&=& A(t) \frac{1}{\sqrt{n}} \sum \limits_{i=1}^n \left\{\frac{\partial \varphi}{\partial \Gamma}\left(-1, \widehat{\Gamma} \right) \right\} \mathcal{D}^{-1} \mathcal{C} \mathbf{B}_i(\mathcal{Z}) + o_p(1),
 \end{eqnarray}
where the third term is due to (\ref{pf-3-8}) and $A(t)$ is the convergent function of $\frac{\partial \widehat{H}(t;b,\zeta,\theta_0)}{\partial \theta}$.

By (\ref{Sec-2-1-3}) and the fact that $\widehat{H}(t;b,\zeta,\theta)$ is a solution of (\ref{EE-SIMEX-1}), we have
\begin{eqnarray*}
&& \sum \limits_{i=1}^n dM_i(t) \\
&=& \sum \limits_{i=1}^n \left[ dN_i(t) -  R_i(t) d \log \left\{  G(\Lambda_\epsilon \left\{H_0(t) + W_i^\top(b,\zeta) \beta_0 \right\} - Z_i^\top \gamma_0) \right\} \right] \\
&=& \sum \limits_{i=1}^n R_i(t) d \log \left\{  G(\Lambda_\epsilon \left\{\widehat{H}(t;b,\zeta,\theta_0) + W_i^\top(b,\zeta) \beta_0 \right\} - Z_i^\top \gamma_0) \right\}  \\
&\ \ &  - \sum \limits_{i=1}^n  R_i(t) d \log \left\{  G(\Lambda_\epsilon \left\{H_0(t) + W_i^\top(b,\zeta) \beta_0 \right\} - Z_i^\top \gamma_0) \right\},
\end{eqnarray*}
which yields that
\begin{eqnarray} \label{pf-4-2}
&& \sqrt{n} \left\{  \widehat{H}(t;b,\zeta,\theta_0) - H_0(t) \right\} \nonumber \\
&=& \frac{1}{\sqrt{n}} \sum \limits_{i=1}^n \int_0^t \frac{B(s,t;W_i(b,\zeta))}{\Psi_i(s;W_i(b,\zeta),H_0,\theta_0)} dM_i(s) + o_p(1).
\end{eqnarray}
Then combining (\ref{pf-4-1}) and (\ref{pf-4-2}) gives
\begin{eqnarray} \label{pf-4-3}
\sqrt{n} \left\{  \widehat{H}(t;b,\zeta,\widehat{\theta}_{SIMEX}) - H_0(t) \right\} = \frac{1}{\sqrt{n}} \sum \limits_{i=1}^n \mathcal{T}_i(t;b,\zeta) + o_p(1),
\end{eqnarray}
where $\mathcal{T}_i(t;b,\zeta) = A(t)  \left\{\frac{\partial \varphi}{\partial \Gamma}\left(-1, \widehat{\Gamma} \right) \right\} \mathcal{D}^{-1} \mathcal{C} \mathbf{B}_i(\mathcal{Z}) + \int_0^t \frac{B(s,t;W_i(b,\zeta))}{\Psi_i(s;W_i(b,\zeta),H_0,\theta_0)} dM_i(s)$. Taking average on both sides of (\ref{pf-4-3}) with respect to $b$ yields
\begin{eqnarray} \label{pf-4-4}
\sqrt{n} \left\{  \widehat{H}(t;\zeta,\widehat{\theta}_{SIMEX}) - H_0(t) \right\} = \frac{1}{\sqrt{n}} \sum \limits_{i=1}^n \mathcal{T}_i(t;\zeta) + o_p(1),
\end{eqnarray}
where $\mathcal{T}_i(t;\zeta) = \frac{1}{B} \sum \limits_{b=1}^B \mathcal{T}_i(t;b,\zeta)$.

Suppose that $\varphi_H(\zeta,\Gamma_H(t))$ is a function with the same conditions in (C7), and $\Gamma_H(t)$ is the associated parameter depending on time $t$. For $t \in [0,\tau]$ and $\zeta \in \mathcal{Z}$, we fit a regression model on $\mathcal{T}_i(t;\zeta)$ and $\varphi_H(\zeta,\Gamma_H(t))$, and derive the estimator of $\Gamma_H(t)$, which is denoted by $\widehat{\Gamma}_H(t)$. Furthermore, similar to the derivations in (\ref{pf-3-7}), we have  
\begin{eqnarray} \label{pf-4-5}
\varphi_H\left( \mathcal{Z}, \widehat{\Gamma}_H(t) \right) - \varphi_H\left( \mathcal{Z}, \Gamma_H(t) \right) \approx \frac{\partial \varphi_H\left( \mathcal{Z}, \Gamma_H(t) \right)}{\partial \Gamma_H(t)} \left\{ \widehat{\Gamma}_H(t) - \Gamma_H(t) \right\}. 
\end{eqnarray}
Let $\mathcal{U}(t) = \frac{\partial \varphi_H\left( \mathcal{Z}, \Gamma_H(t) \right)}{\partial \Gamma_H(t)}$ and $\mathcal{V}(t) = \left\{\frac{\partial \varphi_H\left( \mathcal{Z}, \Gamma_H(t) \right)}{\partial \Gamma_H(t)}\right\}^\top \frac{\partial \varphi\left( \mathcal{Z}, \Gamma_H(t) \right)}{\partial \Gamma_H(t)}$. Combining (\ref{pf-4-4}) and (\ref{pf-4-5}) yields
\begin{eqnarray} \label{pf-4-6}
\sqrt{n} \left\{ \widehat{\Gamma}_H(t) - \Gamma_H(t) \right\} = \frac{1}{\sqrt{n}} \sum \limits_{i=1}^n \mathcal{V}^{-1}(t) \mathcal{U}(t) \mathcal{T}_i(t;\mathcal{Z}) + o_p(1),
\end{eqnarray}
and since the estimator $\widehat{H}_{SIMEX}(t)$ is a predicted value of $\varphi_H(\zeta,\widehat{\Gamma}_H(t))$ by taking $\zeta \rightarrow -1$, then by (\ref{pf-4-6}), we obtain
\begin{eqnarray} \label{pf-4-7}
& & \sqrt{n} \left\{ \widehat{H}_{SIMEX}(t) - H_0(t) \right\} \nonumber \\
&=& \frac{1}{\sqrt{n}} \sum \limits_{i=1}^n \left\{ \frac{\partial \varphi_H(-1,\widehat{\Gamma}_H(t))}{\partial \Gamma_H(t)} \right\} \mathcal{V}^{-1}(t) \mathcal{U}(t) \mathcal{T}_i(t;-1) + o_p(1) \nonumber \\
&\triangleq& \frac{1}{\sqrt{n}} \sum \limits_{i=1}^n \mathcal{H}_i(t) + o_p(1).
\end{eqnarray}
Finally, by the Central Limit Theorem, we conclude that $\sqrt{n} \left\{ \widehat{H}_{SIMEX}(t) - H_0(t) \right\}$ converges to the Gaussian process with mean zero and covariance function $E\left\{ \mathcal{H}_i(t) \mathcal{H}_i(s) \right\}$.
 $\hfill \square$


\section*{References}

\begin{thebibliography}{16}

\bibitem[Anderson et al. (1993)]{Anderson1993} 
Anderson, P.K., Borgan, O., Gill, R.D. and Keiding, N. (1993) {\em Statistical Models Based on Counting Processes}. Springer-Verlag New York.

\bibitem[Bertrand et al. (2017a)]{Bertrand2017a}
Bertrand, A., Legrand, C., Carroll, R. J. Meester, C. D. and Keilegom, I. V. (2017) Inference in a survival cure model with mismeasured covariates using a simulation-extrapolation approach. {\em Biometrika}, 104, 31-50.


\bibitem[Carroll et al. (2006)]{Carroll2006}
Carroll, R. J., Ruppert, D., Stefanski, L. A., and Crainiceanu, C. M. (2006) {\em Measurement Error in Nonlinear Model}. Chapman \& Hall/CRC, New York.



\bibitem[Chen et al. (2017)]{Chen2017}
Chen, C.-M., Shen, P.-S., Wei, J. C.-C. and Lin, L. (2017) A semiparametric mixture cure survival model for left-truncated
and right-censored data. {\em Biometrical Journal}, 59, 270-290.



\bibitem[Cook and Stefaski (1994)]{Cook1994}
Cook, J. R. and Stefaski, L. A. (1994). Simulation-extrapolation estimation in parametric
measurement error models. \textit{Journal of the American Statistical Association}, \textbf{89}, 1314-1328.

\bibitem[Lu and Ying (2004)]{Lu2004}
Lu, W. and Ying, Z. (2004). On semiparametric transformation cure models. {\em Biometrika}, 91, 331–343.


\bibitem[Ma and Yin (2008)]{Ma2008}
Ma, Y. and Yin, G. (2008). Cure rate model with mismeasured covariates under transformation. {\em Journal of the American Statistical Association}, 103, 743–56.



\bibitem[van der Vaart (1998).]{van-der1998}
van der Vaart, A. W. (1998). {\em Asymptotic Statistics}. Cambridge University Press, New York.



















\end{thebibliography}
\end{document}